\documentclass[11pt]{article}
\usepackage{amsmath,latexsym,epsf,mathabx}
\usepackage{epsfig,graphicx,picins,picinpar,subfigure}
\usepackage{pstricks}
\usepackage{fancyvrb}
\usepackage[all]{xy}

\usepackage{ulem}
\begin{document}

\newcommand{\nc}{\newcommand}
\def\PP#1#2#3{{\mathrm{Pres}}^{#1}_{#2}{#3}\setcounter{equation}{0}}
\def\ns{$n$-star}\setcounter{equation}{0}
\def\nt{$n$-tilting}\setcounter{equation}{0}
\def\Ht#1#2#3{{{\mathrm{Hom}}_{#1}({#2},{#3})}\setcounter{equation}{0}}
\def\qp#1{{${(#1)}$-quasi-projective}\setcounter{equation}{0}}
\def\mr#1{{{\mathrm{#1}}}\setcounter{equation}{0}}
\def\mc#1{{{\mathcal{#1}}}\setcounter{equation}{0}}
\def\HD{\mr{Hom}_{\mc{D}}}
\def\HC{\mr{Hom}_{\mc{C}}}
\def\AdT{\mr{Add}_{\mc{T}}}
\def\adT{\mr{add}_{\mc{T}}}
\def\Kb{\mc{K}^b(\mr{Proj}R)}
\def\kb{\mc{K}^b(\mc{P}_R)}
\def\AdpC{\mr{Adp}_{\mc{C}}}

\newtheorem{Th}{Theorem}[section]
\newtheorem{Def}[Th]{Definition}
\newtheorem{Lem}[Th]{Lemma}
\newtheorem{Pro}[Th]{Proposition}
\newtheorem{ex}[Th]{Example}
\newtheorem{Cor}[Th]{Corollary}
\newtheorem{Rem}[Th]{Remark}
\newtheorem{Not}[Th]{Notation}
\newtheorem{Sc}[Th]{}
\def\Pf#1{{\noindent\bf Proof}.\setcounter{equation}{0}}
\def\>#1{{ $\Rightarrow$ }\setcounter{equation}{0}}
\def\<>#1{{ $\Leftrightarrow$ }\setcounter{equation}{0}}
\def\bskip#1{{ \vskip 20pt }\setcounter{equation}{0}}
\def\sskip#1{{ \vskip 5pt }\setcounter{equation}{0}}
\def\mskip#1{{ \vskip 10pt }\setcounter{equation}{0}}
\def\bg#1{\begin{#1}\setcounter{equation}{0}}
\def\ed#1{\end{#1}\setcounter{equation}{0}}
\def\KET{T^{^F\bot}\setcounter{equation}{0}}
\def\KEC{C^{\bot}\setcounter{equation}{0}}

\renewcommand{\thefootnote}{\fnsymbol{footnote}}
\setcounter{footnote}{0}
%
%


\title{\bf  Relative AR-correspondence, co-t-structure and silting pair
\thanks{Supported by the National Natural Science Foundation of China (Grants No.11801004) and the Startup Foundation for Introducing Talent of AHPU (Grant No.2017YQQ016 and 2020YQQ067). }
}
\smallskip
\author{ Peiyu Zhang and Dajun Liu\\ 
\footnotesize ~E-mail:~zhangpy@ahpu.edu.cn, ldjnnu2017004@163.com; \\
\footnotesize School of Mathematics and Physics,  Anhui Polytechnic University, Wuhu, China.  \\
\\ Jiaqun Wei\\ 
\footnotesize ~E-mail:~weijiaqun@njnu.edu.cn;\\
\footnotesize  School of Mathematics Science, Nanjing Normal University, Nanjing, China.
}

\date{}
\maketitle
\baselineskip 15pt
%
%
\begin{abstract}
\vskip 10pt%

As a generalization of tilting pair, which was introduced by Miyashita in \cite{YM}, the notion of silting pair
is introduced in this paper. The authors extends a characterization of tilting modules given by Bazzoni \cite[Theorem~3.11]{BS} to silting pairs,
and proves that there is an one-to-one correspondence between equivalent classes of silting pairs and certain subcategories
which satisfy some conditions. Furthermore, the authors also gives a bijection between equivalent class of silting pairs and bounded above co-t-structure.

\mskip\




\end{abstract}
%
\vskip 30pt

\section{Introduction and Preliminaries}

The tilting theory is well known and plays an important role in the representation theory of
Artin algebra. Beginning with Miyashita \cite{YM1} and in turn Colby and Fuller \cite{CFT}, etc., finitely generated
tilting modules over arbitrary rings were studied. In the study of tilting modules over Artin algebra, Miyashita introduced
the notion of tilting pair in \cite{YM}, i.e., if $A$ is an Artin algebra and $C,~T\in$ $A$-mod (that is, the category of all
finitely generated left $A$-modules), then the pair $(C,~T)$ is tilting if both $C$ and $T$ are selforthogonal such that $T \in \widehat{\mathrm{Add}C}$
and $C\in \widecheck{\mathrm{Add}T}$.

As a generation of tilting theory, silting theory is a hot topic recently.
Silting complexes (i.e., semi-tilting complexes in \cite{WEI}) were first introduced by Keller and Vossieck \cite{BKDV} in order to study t-structures
in the bounded derived category of representations of Dynkin quivers. In \cite{AMV}, the authors introduced the notion of silting module and they proved
that there is a bijection between silting complexes and certain t-structure and co-t-structure in the derived module category.
For articles on the silting theory, interested readers can refer to \cite{A,A1,AMV1,AMV2,BZ,MV,NS,POP,PV,WZ2} etc.
The authors (Zhang and Wei) introduced three notions of cosilting complexes, cosilting modules (simultaneously and independently with \cite{BF})
and AIR-cotilting modules, and it is proved that they are the same in \cite{WZ1}.

In this paper the authors introduces the new notion of silting pairs, gives a characterization of silting pairs, proves that
there is an one-to-one correspondence between equivalent classes of silting pairs and certain subcategories
which satisfies some conditions, and obtains a bijection between silting pair and bounded above co-t-structure.

Let $R$ be an ring. We denote by $R$-Mod the category of all left $R$-modules.
 We denote by $\mathrm{Proj}R$ the class of all projective $R$-modules. The notation
$K^{b}(\mathrm{Proj}R)$ denotes the homotopy category of bounded complexes of projective modules.
Denote by $\mathcal{D}(R)$ (resp. $\mathcal{D}^{b}(R)$) the category of the unbounded (resp. bounded) derived category of complexes
in $R$-Mod.
Throughout this paper we consider a triangulated category $\mathcal{D}$ with [1] the shift functor.
Assume that $\mathcal{B}$ is a full subcategory of $\mathcal{D}$. Recall that $\mathcal{B}$ is closed under extensions if for any triangle $X\to Y\to Z\to$ in  $\mathcal{D}$ with
$X, Z \in \mathcal{B}$, we have $Y\in \mathcal{B}$. The subcategory
$\mathcal{B}$ is suspended (resp. cosuspended) if it is closed under extensions and under the functor $[1]$ (resp. $[-1]$).
An object $M\in \mathcal{D}$ has a $\mathcal{B}$-resolution (resp., $\mathcal{B}$-coresolution) with the length
at most $m$ ($m\ge 0$), if there are
triangles $M_{i+1}\to X_i\to M_i\to$ (resp., $M_i\to X_i\to
M_{i+1}\to$) with $0\le i\le m$ such that $M_0=M$, $M_{m+1}=0$
and each $X_i\in \mathcal{B}$. In this case, we denoted by $\mathcal{B}$-res.dim$(M)\le m$
(resp., $\mathcal{B}$-cores.dim$(M)\le m$).

Next, given a subcategory $\mathcal{B}$ of $\mathcal{D}$,  we give the following notations which are widely used in the tilting theory,
where $n\ge 0$ and $m$ is an integer.

\bg{verse}

$(\widehat{\mathcal{B}})_n=\{L\in \mathcal{D} |
\mathcal{B}$-res.dim$(L)\le n\}$.

$(\widecheck{\mathcal{B}})_n=\{L\in\mathcal{D}\ |\
\mathcal{B}$-cores.dim$(L)\le n\}$.

$\widehat{\mathcal{B}}\ \ \ \ =\{L\in\mathcal{D}\ |\ L\in
(\widehat{\mathcal{B}})_n$ for some $n\}$.

$\widecheck{\mathcal{B}}\ \ \ \ =\{L\in\mathcal{D}\ |\ L\in
(\widecheck{\mathcal{B}})_n$ for some $n\}$.

${\mathcal{B}}^{\bot_{>m}}=\{N\in \mathcal{D}\ |\
\mathrm{Hom}_{\mathcal{D}}(M,N[i])=0$ for all $M\in \mathcal{B}$ and all $i>m\}$.

$^{\bot_{>m}}{\mathcal{B}}=\{N\in \mathcal{D}\ |\
\mathrm{Hom}_{\mathcal{D}}(N,M[i])=0$ for all $M\in \mathcal{B}$ and all $i>m\}$.

$^{\bot_{m}}{\mathcal{B}}=\{N\in \mathcal{D}\ |\
\mathrm{Hom}_{\mathcal{D}}(N,M[m])=0$ for all $M\in \mathcal{B}$ $\}$.

${\mathcal{B}}^{\bot_{\gg 0}}\ =\{N\in \mathcal{D}\ |\
N\in{\mathcal{B}}^{\bot_{>m}}$ for some $m\}$.


\ed{verse}

\noindent Moreover, for an object $M\in \mathcal{D}$, we define the class $\mathrm{Add}M=\{ N\in \mathcal{D}|$ there exists $L$ such that $N\oplus L=M^{(X)}$ for some set $X$ $\}$.
If the category $\mathcal{B}=\mathrm{Add}M$, then we simply replace the notation ${\mathcal{B}}^{\bot_{>m}}$
with ${M}^{\bot_{>m}}$ (the other situations are similar). Also, we need the following two notions.

\bg{verse}

${\mc{X}_{M}}\ \ =\{N\in {^{\bot_{>0}}}M\ |$ there are
triangles $N_{i}\to M_i\to N_{i+1}\to$  such that
$N_0=N$, $N_i\in {^{\bot_{>0}}M}$ and $M_i\in\mathrm{Add}M$ for
all $i\ge 0\}$.

${_{M}\mc{X}}\ =\{N\in {M}^{\bot_{>0}}\ |$ there are
triangles $N_{i+1}\to M_i\to N_i\to$  such that $N_0=N$,
$N_i\in {M}^{\bot_{>0}}$ and $M_i\in \mathrm{Add}M$ for all $i\ge 0\}$.

\ed{verse}
%

\section{Silting pairs}

Recall that an object $M$ is said to be semi-selforthogonal if $\mathrm{Add}M\subseteq M^{\bot_{>0}}$.
Firstly, we collect some important results about semi-selforthogonal objects from the second section of \cite{WEI} in the following
lemma.

\bg{Lem}\label{WEI}
Let $M$ be semi-selforthogonal in $\mathcal{D}$. Then the following results hold.

$(1)$ $\mathcal{X}_{M}$ is closed under extensions, direct summands and $[-1]$;

$(1')$ $_{M}\mathcal{X}$ is closed under extensions, direct summands and $[1]$;

$(2)$ $\mathrm{Hom}_{\mathcal{D}}(\widecheck{\mathrm{Add}M}, ~M^{\bot_{>0}}[i>0])=0$;

$(2')$ $\mathrm{Hom}_{\mathcal{D}}(^{\bot_{>0}}M, ~\widehat{\mathrm{Add}M}[i>0])=0$;

$(3)$ $(\widehat{\mathrm{Add}M})_{n}=$ $_{M}\mathcal{X}\bigcap ^{\bot_{>n}}(_{M}\mathcal{X})
=~_{M}\mathcal{X}\bigcap ^{\bot_{>n}}(M^{\bot_{>0}})$;

$(3')$ $(\widecheck{\mathrm{Add}M})_{n}=\mathcal{X}_{M}\bigcap (\mathcal{X}_{M})^{\bot_{>n}}
=\mathcal{X}_{M}\bigcap (^{\bot_{>0}}M)^{\bot_{>n}}$;

$(4)$ Both $\widecheck{\mathrm{Add}M}$ and $\widehat{\mathrm{Add}M}$ are closed under extensions and direct summands.

$(5)$ $\mathrm{Add}M=\widehat{\mathrm{Add}M}\bigcap$ $^{\bot_{>0}}M=\widecheck{\mathrm{Add}M}\bigcap$ $M^{\bot_{>0}}$
\ed{Lem}

We also need the following results.

\bg{Lem}\label{LBST}
Let $M$ be semi-selforthogonal in $\mathcal{D}$. Given any triangle $X\longrightarrow Y\longrightarrow Z\longrightarrow$, then the following results hold.

$(1)$ If $Y, ~Z \in$ $_{M}\mathcal{X}$
and $\mathrm{Hom}_{\mathcal{D}}(M, X[1])=0$, then $X \in$ $_{M}\mathcal{X}$;

$(2)$ If $X, ~Y \in$ $\mathcal{X}_{M}$
and $\mathrm{Hom}_{\mathcal{D}}(Z, M[1])=0$, then $Z \in$ $\mathcal{X}_{M}$.
\ed{Lem}

\Pf. (1) Since $Z \in$ $_{M}\mathcal{X}$, there is a triangle
$Z'\longrightarrow M'\longrightarrow Z\longrightarrow$ with $M'\in \mathrm{Add}M$ and $Z' \in$ $_{M}\mathcal{X}$.
Thus we have the following triangles commutative diagram:

$$\xymatrix{
&Z'\ar[d]\ar@^{=}[r]&Z'\ar[d]\\
X\ar[r]\ar@^{=}[d]&U\ar[r]\ar[d]&M'\ar[d]\ar[r]&\\
X\ar[r]&Y\ar[r]\ar[d]&Z\ar[r]\ar[d]&\\
&&
}$$
By the lemma \ref{WEI} (1'), we have that $U\in$ $_{M}\mathcal{X}$ from the second column in diagram above.
Since $\mathrm{Hom}_{\mathcal{D}}(M, X[1])=0$, the second row is split.
It follows from the lemma \ref{WEI} (1') that $X\in$ $_{M}\mathcal{X}$.

(2) The proof is dual to (1).
\ \hfill $\Box$

\vspace{6pt}

Let $\mathcal{W}\subseteq \mathcal{D}$ be a subcategory closed under finite direct sums and
summands. We say that $M$ is a relative generator of $\mathcal{W}$, if $M \in \mathcal{W}$ and for
any $W\in \mathcal{W}$, there is a triangle $W'\longrightarrow M_{W}\longrightarrow W\longrightarrow$
with $M_{W}\in \mathrm{Add}M$ and $W'\in \mathcal{W}$. Dually, we can define the relative cogenerator.
For example, if $M$ is semi-selforthogonal, then $_{M}\mathcal{X}$ and $\widehat{\mathrm{Add}M}$ (respectively, $\mathcal{X}_{M}$
and $\widecheck{\mathrm{Add}M}$) both have a relative generator (respectively, cogenerator) $M$.

\mskip\
The following lemma comes from the idea in \cite[Lemma~2.2]{WX}.

\bg{Lem}\label{ATWC}
Let $\mathcal{W}$ be a subcategory of $\mathcal{D}$. Assume that $\mathcal{W}$ contains a semi-selforthogonal relative
generator $M$ and is closed under extensions, finite direct sums and direct
summands. If there are triangles $X_{i+1}\longrightarrow Y_{i}\longrightarrow X_{i}\longrightarrow$ with each
$Y_{i}\in \mathcal{W}$ for any $1\leq i\leq n$,
then there exist $U_{n}$ and $V_{n+1}$ satisfying the following two conditions:

$(1)$ there is a triangle $U_{n}\longrightarrow V_{n+1} \longrightarrow X_{n+1}\longrightarrow$
with $U_{n}\in \mathcal{W}$;

$(2)$  there are triangles $V_{i+1}\longrightarrow Z_{i}\longrightarrow V_{i}\longrightarrow$
with $Z_{i} \in \mathrm{Add}M$ for any $1\leq i\leq n$ and $V_{1}=X_{1}$.

Moreover, if $X_{1}\in \mathcal{W}$, then there is a triangle
$U_{n}\longrightarrow V_{n+1} \longrightarrow X_{n+1}\longrightarrow$ with $U_{n}\in \mathcal{W}$ and $V_{n+1}\in (\widecheck{\mathrm{Add}M})_{n}$.
\ed{Lem}

\Pf. By induction on $n$. When $n=1$, we have a triangle $X_{2}\longrightarrow Y_{1}\longrightarrow X_{1}\longrightarrow$
with $Y_{1}\in \mathcal{W}$.
Thus there is a triangle $U_{1}\longrightarrow Z_{1}\longrightarrow Y_{1}\longrightarrow$
with $Z_{1}\in \mathrm{Add}M$ and $U_{1}\in \mathcal{W}$.
So we have the following triangles commutative diagram:

$$\xymatrix{
U_{1}\ar[d]\ar@^{=}[r]&U_{1}\ar[d]\\
V_{2}\ar[r]\ar[d]&Z_{1}\ar[r]\ar[d]&X_{1}\ar@^{=}[d]\ar[r]&\\
X_{2}\ar[r]\ar[d]&Y_{1}\ar[r]\ar[d]&X_{1}\ar[r]&\\
&&
}$$
It is easy to see that the left column and the middle row are just the desired triangles.

Assume that the result holds for $n-1$. Then there exist $U_{n-1}$ and $V_{n}$ such that there is a triangle
$U_{n-1}\longrightarrow V_{n} \longrightarrow X_{n}\longrightarrow$ with $U_{n-1}\in \mathcal{W}$
and $V_{i+1}\longrightarrow Z_{i}\longrightarrow V_{i}\longrightarrow$
with $Z_{i} \in \mathrm{Add}M$ for any $1\leq i\leq n-1$ and $V_{1}=X_{1}$.
Consequently, we can obtain the following triangles commutative diagram:

$$\xymatrix{
&U_{n-1}\ar[d]\ar@^{=}[r]&U_{n-1}\ar[d]\\
X_{n+1}\ar[r]\ar@^{=}[d]&A\ar[r]\ar[d]&V_{n}\ar[d]\ar[r]&\\
X_{n+1}\ar[r]&Y_{n}\ar[r]\ar[d]&X_{n}\ar[r]\ar[d]&\\
&&
}$$
It follows, since $\mathcal{W}$ is closed under extensions, that $A$ is in $\mathcal{W}$.
So we have a triangle $U_{n}\longrightarrow Z_{n}\longrightarrow A\longrightarrow$ with
$Z_{n}\in \mathrm{Add}M$ and $U_{n}\in \mathcal{W}$.
Then we can construct the following triangles commutative diagram:

$$\xymatrix{
U_{n}\ar[d]\ar@^{=}[r]&U_{n}\ar[d]\\
V_{n+1}\ar[r]\ar[d]&Z_{n}\ar[r]\ar[d]&V_{n}\ar@^{=}[d]\ar[r]&\\
X_{n+1}\ar[r]\ar[d]&A\ar[r]\ar[d]&V_{n}\ar[r]&\\
&&
}$$
So we complete this proof.

Moreover, if $X_{1}\in \mathcal{W}$,
then there is a triangle $X_{1}'\longrightarrow M'\longrightarrow X_{1}\longrightarrow$
with $M'\in \mathrm{Add}M$ and $X_{1}'\in \mathcal{W}$.
We have the following triangles commutative diagram:

$$\xymatrix{
&X_{1}'\ar[d]\ar@^{=}[r]&X_{1}'\ar[d]\\
X_{2}\ar[r]\ar@^{=}[d]&B\ar[r]\ar[d]&M'\ar[d]\ar[r]&\\
X_{2}\ar[r]&Y_{1}\ar[r]\ar[d]&X_{1}\ar[r]\ar[d]&\\
&&
}$$
It follows, since $\mathcal{W}$ is closed under extensions, that $B\in \mathcal{W}$.
Replace the triangle $X_{2}\longrightarrow Y_{1}\longrightarrow X_{1}\longrightarrow$ with
the triangle $X_{2}\longrightarrow B\longrightarrow M'\longrightarrow$, we can
obtain the conclusion from (1) and (2).
\ \hfill $\Box$

\bg{Lem}\label{ATWI}
Assume that M is a semi-selforthogonal object in $\mathcal{D}$.
If there are triangles $X_{i-1}\longrightarrow Y_{i}\longrightarrow X_{i}\longrightarrow$ with each
$Y_{i}\in (\widecheck{\mathrm{Add}M})_{n}$, for any $1\leq i\leq n$, then there exist $U_{n}$ and $V_{n}$ satisfying the following two conditions:

$(1)$ there is a triangle $X_{n}\longrightarrow V_{n} \longrightarrow U_{n}\longrightarrow$ with $U_{n}\in (\widecheck{\mathrm{Add}M})_{n-1}$;

$(2)$  there are triangles $V_{i-1}\longrightarrow Z_{i}\longrightarrow V_{i}\longrightarrow$
with $Z_{i} \in \mathrm{Add}M$ for any $1\leq i\leq n$ and $V_{0}=X_{0}$.

\ed{Lem}

\Pf. The proof is dual to Lemma \ref{ATWC}.
\ \hfill $\Box$

%
%
%

\bg{Pro}\label{add}
Let $M$ be a semi-selforthogonal object in $\mathcal{D}$. Then $\widecheck{_{M}\mathcal{X}}$ is a
triangulated subcategory of $\mathcal{D}$.

\ed{Pro}

\Pf. By Lemma 2.1 (2) in \cite{WEI1} and Lemma 2.1 (1').
\ \hfill $\Box$

\mskip\

Next, we give the definition of silting pair and some useful properties.

\bg{Def}\label{}
A pair $(C,~S)$ is silting if it satisfies the following conditions:

$(1)$ $C$ is semi-selforthogonal;

$(2)$ $S$ is semi-selforthogonal;

$(3)$ $C\in \widecheck{\mathrm{Add}S}$;

$(4)$ $S\in \widehat{\mathrm{Add}C}$.

\ed{Def}

Let $R$ be a ring and $\mathcal{D}^{b}(R)$ be its bounded derived category.
We say that a complex $T$ is said to be small silting \cite[Definition~3.1]{WEI} if it satisfies the following conditions:
$(i)$ $T\in \mathcal{K}^b(\mathrm{proj}R)$, $(ii)$ $T$ is semi-selforthogonal, and $(iii)$ $\mathcal{K}^b(\mathrm{proj}R) = \langle \mathrm{add}T\rangle$,
i.e., $\mathcal{K}^b(\mathrm{proj}R)$ coincides with the  smallest triangulated subcategory containing $\mathrm{add}T$, where
$\mathrm{add}T=\{ M\in \mathcal{D}|$ there exists $N$ such that $N\oplus M=T^{(X)}$ for some finite set $X$ $\}$.
From the definition of silting pair, the following are obvious:

$(1)$ If $C$ is semi-selforthogonal, then $(C, C)$ is a silting pair;

$(2)$ If $T\in\mathcal{D}^{\leq0}$ (the complexes whose homologies are concentrated on non-positive terms)
is a small silting complex, then $(R, T)$ is a silting pair by the small version of Theorem 3.5 in \cite{WEI};

$(3)$ If $T$ is a tilting complex, then $(R, T)$ is a silting pair;

$(4)$ If $M$ is a tilting R-module, then $(R, M)$ is a silting pair by the small version of Corollary 3.7 in \cite{WEI}.

\bg{ex}\label{quiver}
Let $A$ be the path $K$-algebra given by the linear quiver
$\xymatrix{ 1\ar[r] & 2 \ar[r] & 3 },$ where $K$ is an algebraically closed field. Then its AR quiver is given by
$$
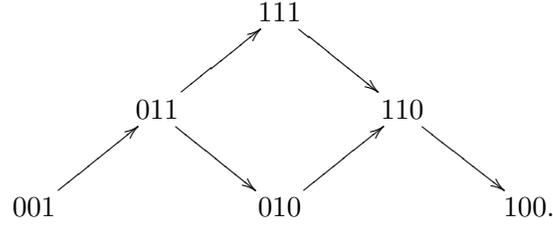
\xymatrix{
&&111\ar[dr]&&\\
&011\ar[ur]\ar[dr]&&110\ar[dr]&\\
001\ar[ur] &&010\ar[ur]&&100.
}$$
In $\mathcal{D}^{b}(A)$, we take $C=110\oplus100$, $S=100\oplus010[1]$,
it is easy to see that both $C$ and $S$ are semi-selforthogonal, and there is a triangle
$110\longrightarrow 100\longrightarrow 010[1] \longrightarrow$. So the pair $(C,~S)$ is silting.

\ed{ex}

\noindent
From the following the lemma, we will give the definition of $n$-silting pair.

\bg{Lem}\label{twon}
Let the pair $(C,~S)$ be silting and m be a positive integer. Then $C\in (\widecheck{\mathrm{Add}S})_{m}$
if and only if $S\in (\widehat{\mathrm{Add}C})_{m}$.

\ed{Lem}

\Pf. $(\Longleftarrow)$: Assume that $C\in (\widecheck{\mathrm{Add}S})_{n}$ and $n>m$.
Thus there are triangles $C_{i}\longrightarrow S_{i}\longrightarrow C_{i+1}\longrightarrow$
with $S_{i}\in \mathrm{Add}S$ for any $0\leq i\leq n-1$, where $C_{0}=C$ and $C_{n}=S_{n}$.
Using the functor $\mathrm{Hom}_{\mathcal{D}}(C_{n},~-)$ to these triangles, we have that
$\mathrm{Hom}_{\mathcal{D}}(C_{n},C_{n-1}[1])\cong \mathrm{Hom}_{\mathcal{D}}(C_{n},C_{0}[n])$ by dimension shifting.
By the lemma \ref{WEI} (4) and (3), we have that
$\mathrm{Hom}_{\mathcal{D}}(C_{n},C_{n-1}[1])\cong \mathrm{Hom}_{\mathcal{D}}(C_{n},C_{0}[n])=0$.
Thus the triangle $C_{n-1}\longrightarrow S_{n-1}\longrightarrow C_{n}\longrightarrow$ is split.
It follows that $C\in (\widecheck{\mathrm{Add}S})_{n-1}$.
Repeating this process, we finally obtain that $C\in (\widecheck{\mathrm{Add}S})_{m}$.

$(\Longrightarrow)$: The proof is dual to the above.
\ \hfill $\Box$

\bg{Def}\label{n-silting}
A pair $(C,~S)$ is n-silting if it satisfies the following two conditions:

$(1)$ the pair $(C,~S)$ is silting;

$(2)$ $C\in (\widecheck{\mathrm{Add}S})_{n}$.
\ed{Def}

From Lemma \ref{twon}, it is easy to see that we can replace the condition (2) in Definition \ref{n-silting}
with the condition $S\in (\widehat{\mathrm{Add}C})_{n}$.

\bg{Lem}\label{LIIL}
Let $S$ and C be semi-selforthogonal in $\mathcal{D}$.

$(1)$ If $C\in \widecheck{\mathrm{Add}S}$, then $S^{\bot_{>0}}\subseteq C^{\bot_{>0}}$;

$(2)$ If $S\in$ $_{C}\mathcal{X}$ and $C\in \widecheck{\mathrm{Add}S}$, then
$_{S}\mathcal{X}\subseteq$ $_{C}\mathcal{X}$.

\ed{Lem}

\Pf. (1) Easily.

(2) Taking $T\in$ $_{S}\mathcal{X}$, there are triangles
$T_{i+1}\longrightarrow S_{i}\longrightarrow T_{i}\longrightarrow$
with $S_{i}\in \mathrm{Add}S$ for any $i\geq0$, where $T_{0}=T$ and $T_{j}\in$ $_{S}\mathcal{X}$ for any $j>0$.
Since $_{C}\mathcal{X}$ is closed under extensions, $S_{i}\in$ $_{C}\mathcal{X}$ for any $i\geq1$.
Note that $_{C}\mathcal{X}$ is closed under finite direct sums by Lemma 3.11 from \cite{WEI}.
By the lemma \ref{ATWC}, thus we can obtain a triangle $U_{0}\longrightarrow V_{1}\longrightarrow T_{1}\longrightarrow$
with $U_{0}\in$ $_{C}\mathcal{X}$ such that there is a triangle
$V_{1}\longrightarrow C_{0}\longrightarrow T_{0}\longrightarrow$ with $C_{0}\in \mathrm{Add}C$.
Thus, we have the following triangles commutative diagram:

$$\xymatrix{
&T_{2}\ar[d]\ar@^{=}[r]&T_{2}\ar[d]\\
U_{0}\ar[r]\ar@^{=}[d]&Y\ar[r]\ar[d]&S_{1}\ar[d]\ar[r]&\\
U_{0}\ar[r]&V_{1}\ar[r]\ar[d]&T_{1}\ar[r]\ar[d]&\\
&&
}$$
Since $S_{i}$ and $U_{0}$ are in $_{C}\mathcal{X}$, $Y\in$ $_{C}\mathcal{X}$ from the second in above diagram.
Then there is a triangle
$U_{1}\longrightarrow C_{1}\longrightarrow Y\longrightarrow$ with $C_{1}\in \mathrm{Add}C$ and $U_{1}\in$ $_{C}\mathcal{X}$.
So we have the following triangle commutative diagram:

$$\xymatrix{
U_{1}\ar[d]\ar@^{=}[r]&U_{1}\ar[d]&\\
V_{2}\ar[r]\ar[d]&C_{1}\ar[r]\ar[d]&V_{1}\ar@^{=}[d]\ar[r]&\\
T_{2}\ar[d]\ar[r]&Y\ar[r]\ar[d]&V_{1}\ar[r]&\\
&&
}$$
Now, we have obtained two triangles
$V_{1}\longrightarrow C_{0}\longrightarrow T_{0}\longrightarrow$ and
$V_{2}\longrightarrow C_{1}\longrightarrow V_{1}\longrightarrow$.
Repeating the process for the triangle $U_{1}\longrightarrow V_{2}\longrightarrow T_{2}\longrightarrow$,
we can get a triangle $V_{3}\longrightarrow C_{2}\longrightarrow V_{2}\longrightarrow$.
Keep repeating, we can know that $T\in$ $_{C}\mathcal{X}$. i.e., $_{S}\mathcal{X}\subseteq$ $_{C}\mathcal{X}$.
\ \hfill $\Box$

\mskip\

In order to give a characterization of $n$-silting pairs, we need to introduce the following subcategory.
Let $\mathcal{B}$ and $\mathcal{C}$ be two subcategories of $\mathcal{D}$ and let $n$ be a positive integer. We denote by

\bg{verse}

$\mathrm{Pres}^{n}_{\mathcal{C}}(\mathcal{B})$ := $\{X\in\mathcal{D}|$ there exist triangles
$X_{i}\longrightarrow T_{i-1}\longrightarrow X_{i-1}\longrightarrow$
with $T_{i}\in \mathcal{B}$ for any $1\leq i\leq n$, where $X_{0}=X$ and $X_{n}\in \mathcal{C} \}$.

\ed{verse}

\noindent If the category $\mathcal{B}=\mathrm{Add}T$, then we simply replace the notation $\mathrm{Pres}^{n}_{\mathcal{C}}(\mathcal{B})$
with $\mathrm{Pres}^{n}_{\mathcal{C}}(T)$.

\bg{Pro}\label{three}
Assume that the pair $(C,~S)$ is $n$-silting and $T\in \mathcal{D}$. Then the following statements are equivalent.

$(1)$ $T\in$ $_{S}\mathcal{X}$;

$(2)$ $T\in S^{\bot_{>0}}\bigcap$ $_{C}\mathcal{X}$;

$(3)$ $T\in \mathrm{Pres}^{n}_{_{C}\mathcal{X}}(S)$.

\ed{Pro}

\Pf. $(1)\Longrightarrow(2)$ By lemma \ref{LIIL} (2) we only need to prove that $S\in$ $_{C}\mathcal{X}$.
It follows, by Lemma \ref{WEI}(3) that $S\in$ $_{C}\mathcal{X}$.

$(2)\Longrightarrow(3)$ Assume that $T\in S^{\bot_{>0}}\bigcap$ $_{C}\mathcal{X}$.
Then there is a triangle $T_{1}\longrightarrow S^{(X)}\stackrel{u}{\longrightarrow} T\longrightarrow$
with $u$ evaluation map, where $X=\mathrm{Hom}_{\mathcal{D}}(S,~T)$ and $\mathrm{Hom}_{\mathcal{D}}(S,~u)$ is surjective.
It is easy to prove that $T_{1}\in S^{\bot_{>0}}\subseteq C^{\bot_{>0}}$.
Note that $S\in \widehat{\mathrm{Add}C}\subseteq$ $_{C}\mathcal{X}$ by Lemma \ref{WEI}.
Thus $T_{1}\in $ $_{C}\mathcal{X}$ by Lemma \ref{LBST}.
So $T_{1}\in S^{\bot_{>0}}\bigcap$ $_{C}\mathcal{X}$.
Repeating this process, we have that $T\in \mathrm{Pres}^{n}_{_{C}\mathcal{X}}(S)$.

$(3)\Longrightarrow(2)$ For any $T\in \mathrm{Pres}^{n}_{_{C}\mathcal{X}}(S)$,
there are triangles $M\longrightarrow S_{n}\longrightarrow T_{n-1}\longrightarrow$ and
$T_{i}\longrightarrow S_{i}\longrightarrow T_{i-1}\longrightarrow$
with $M \in$ $_{C}\mathcal{X}$ and $S_{i}\in \mathrm{Add}S$ for any $1\leq i\leq n-1$, where $T_{0}=T$.
Note $\mathrm{Add}S\subseteq$ $_{C}\mathcal{X}$.
We have that $T\in$ $_{C}\mathcal{X}$ by Lemma \ref{WEI} (1').
Applying the functor $\mathrm{Hom}_{\mathcal{D}}(S,~-)$ to these triangles, we have that
$\mathrm{Hom}_{\mathcal{D}}(S,M[i])\cong \mathrm{Hom}_{\mathcal{D}}(S,T[i-n])$ for $i>n$.
Since $S\in \widehat{\mathrm{Add}C}$, there exist triangles
$Q_{i}\longrightarrow C_{i-1}\longrightarrow Q_{i-1}\longrightarrow$
with $C_{i}\in \mathrm{Add}C$ for any $1\leq i\leq n+1$, where $Q_{0}=S$ and $Q_{n+1}=0$.
Applying the functor $\mathrm{Hom}_{\mathcal{D}}(-,~M)$ to these triangles, we have that
$\mathrm{Hom}_{\mathcal{D}}(S,M[i])\cong \mathrm{Hom}_{\mathcal{D}}(C_{n},M[i-n])=0$ for $i>n$.
It follows that $\mathrm{Hom}_{\mathcal{D}}(S,T[i>0])\cong \mathrm{Hom}_{\mathcal{D}}(C_{n},M[i>0])=0$, i.e., $T\in S^{\bot_{>0}}$.
Consequently, (2) holds.

$(2)\Longrightarrow(1)$ In fact, this proof has been completed in $(2)\Longrightarrow(3)$.
\ \hfill $\Box$

\bg{Lem}\label{supplement}
Assume that $C$ is semi-selforthogonal and $S\in \mathcal{D}$.
If $\mathrm{Pres}^{n}_{_{C}\mathcal{X}}(S)=S^{\bot_{>0}}\bigcap$ $_{C}\mathcal{X}$ for some integer $n\geq1$,
then $\mathrm{Add}S\subseteq S^{\bot_{>0}}\bigcap$ $_{C}\mathcal{X}$. Specially, $S$ is semi-selforthogonal.

\ed{Lem}

\Pf. Note that $S^{(X)}\in\mathrm{Pres}^{n}_{_{C}\mathcal{X}}(S)$ for any set $X$. Since both $S^{\bot_{>0}}$ and $_{C}\mathcal{X}$ are closed
under direct summands, $\mathrm{Add}S\subseteq S^{\bot_{>0}}\bigcap$ $_{C}\mathcal{X}$.
\ \hfill $\Box$

\bg{Pro}\label{nn1}
Assume that $C$ is semi-selforthogonal and $S\in \mathcal{D}$.
If $\mathrm{Pres}^{n}_{_{C}\mathcal{X}}(S)=S^{\bot_{>0}}\bigcap$ $_{C}\mathcal{X}$,
then $\mathrm{Pres}^{n}_{_{C}\mathcal{X}}(S)=\mathrm{Pres}^{n+1}_{_{C}\mathcal{X}}(S)$.

\ed{Pro}

\Pf. For any $M\in\mathrm{Pres}^{n+1}_{_{C}\mathcal{X}}(S)$, there are triangles
$M_{i+1}\longrightarrow S_{i+1}\longrightarrow M_{i}\longrightarrow$ with $S_{i+1}\in \mathrm{Add}S$
for any $0\leq i\leq n$, where $M_{0}=M$ and $M_{n+1}\in~_{C}\mathcal{X}$.
Consider the triangle $M_{n+1}\longrightarrow S_{n+1}\longrightarrow M_{n}\longrightarrow$, by Lemma \ref{WEI} (1') and Lemma \ref{supplement},
we can obtain that $M_{n}\in$ $_{C}\mathcal{X}$. i.e., $M\in\mathrm{Pres}^{n}_{_{C}\mathcal{X}}(S)$.

Conversely, for any $T\in \mathrm{Pres}^{n}_{_{C}\mathcal{X}}(S)$, there is a triangle
$T_{1}\longrightarrow S_{1}\longrightarrow T\longrightarrow$ with $S_{1}\in \mathrm{Add}S$ and
$T_{1}\in \mathrm{Pres}^{n-1}_{_{C}\mathcal{X}}(S)$, and there is a triangle
$M\longrightarrow S^{(I)}\stackrel{u}{\longrightarrow} T\longrightarrow$ with $u$ evaluation map.
It is easy to verify that $M\in S^{\bot_{>0}}$.
By the lemma \ref{supplement}, both $S_{1}$ and $S^{(I)}$ are in $_{C}\mathcal{X}$.
We have the following triangles commutative diagram:

$$\xymatrix{
&T_{1}\ar[d]\ar@^{=}[r]&T_{1}\ar[d]&\\
M\ar[r]\ar@^{=}[d]&N\ar[r]\ar[d]&S_{1}\ar[d]\ar[r]&\\
M\ar[r]&S^{(I)}\ar[r]\ar[d]&T\ar[d]\ar[r]&\\
&&
}$$
It follows from $M\in S^{\bot_{>0}}$ that the second row is split.
From the first part of the proof, it is easy to see that $\mathrm{Pres}^{n-1}_{_{C}\mathcal{X}}(S)\subseteq$ $_{C}\mathcal{X}$.
Note that $T_{1}\in \mathrm{Pres}^{n-1}_{_{C}\mathcal{X}}(S)\subseteq$ $_{C}\mathcal{X}$.
From the second column, we have that $N\in$ $_{C}\mathcal{X}$ since $_{C}\mathcal{X}$ is closed under extensions.
Since $_{C}\mathcal{X}$ is closed under direct summands, $M\in$ $_{C}\mathcal{X}$. Consequently,
$M\in$ $_{C}\mathcal{X}\bigcap S^{\bot_{>0}}=$
$\mathrm{Pres}^{n}_{_{C}\mathcal{X}}(S)$.
So we have that $T\in \mathrm{Pres}^{n+1}_{_{C}\mathcal{X}}(S)$.
\ \hfill $\Box$

\bg{Pro}\label{PRES}
Assume that $C$ is semi-selforthogonal and $S\in \mathcal{D}$.
If $\mathrm{Pres}^{n}_{_{C}\mathcal{X}}(S)=S^{\bot_{>0}}\bigcap$ $_{C}\mathcal{X}$,
then the following conclusions hold.

$(1)$ $\mathrm{Pres}^{n}_{_{C}\mathcal{X}}(\mathrm{Pres}^{n}_{_{C}\mathcal{X}}(S))
=\mathrm{Pres}^{n}_{_{C}\mathcal{X}}(S)$;

$(2)$ $C\in (\widecheck{\mathrm{Add}S})_{n}$, and

$(3)$ $S\in (\widehat{\mathrm{Add}C})_{n}$

\ed{Pro}

\Pf. (1) By the lemma \ref{supplement}, we only need to prove that
$\mathrm{Pres}^{n}_{_{C}\mathcal{X}}(\mathrm{Pres}^{n}_{_{C}\mathcal{X}}(S))
\subseteq\mathrm{Pres}^{n}_{_{C}\mathcal{X}}(S)$.
Taking any $T\in \mathrm{Pres}^{n}_{_{C}\mathcal{X}}(\mathrm{Pres}^{n}_{_{C}\mathcal{X}}(S))$,
there are triangles $N\longrightarrow M_{n}\longrightarrow T_{n-1}\longrightarrow$ and
$T_{i}\longrightarrow M_{i}\longrightarrow T_{i-1}\longrightarrow$ with
$M_{i}\in \mathrm{Pres}^{n}_{_{C}\mathcal{X}}(S)$ for any $1\leq i\leq n$ and $N\in$ $_{C}\mathcal{X}$,
where $T_{0}=T$. For any $X \in \mathrm{Pres}^{n}_{_{C}\mathcal{X}}(S)$, by Proposition \ref{nn1}, we have a
$Y\longrightarrow S'\longrightarrow X\longrightarrow$ with $Y \in \mathrm{Pres}^{n}_{_{C}\mathcal{X}}(S)$
and $S'\in \mathrm{Add}S$. i.e.,
$S^{\bot_{>0}}\bigcap$ $_{C}\mathcal{X}=\mathrm{Pres}^{n}_{_{C}\mathcal{X}}(S)$
have a semi-selforthogonal relative generator $S$.
Consequently, $\mathrm{Pres}^{n}_{_{C}\mathcal{X}}(S)=S^{\bot_{>0}}\bigcap$ $_{C}\mathcal{X}$
satisfies the assumption of Lemma \ref{ATWC}. Then there is a triangle
$U\longrightarrow V\longrightarrow N\longrightarrow$ with
$U\in\mathrm{Pres}^{n}_{_{C}\mathcal{X}}(S)$ and there are triangles
$V_{i+1}\longrightarrow Z_{i}\longrightarrow V_{i}\longrightarrow$
with $Z_{i} \in \mathrm{Add}S$ for any $1\leq i\leq n$,
where $V_{n+1}=V$ and $V_{1}=T$. Note that both $U$ and $N$ are in $_{C}\mathcal{X}$, then $V\in$ $_{C}\mathcal{X}$.
So we have that $T\in\mathrm{Pres}^{n}_{_{C}\mathcal{X}}(S)$ by the definition.

(2) Note that $0\in \mathrm{Add}S$ and $C\in$ $_{C}\mathcal{X}$. So $C[n]\in \mathrm{Pres}^{n}_{_{C}\mathcal{X}}(S)$.
i.e., there are some triangles $C[i] \longrightarrow 0\longrightarrow C[i+1]\longrightarrow$ with $0\leq i\leq n-1$.
It is easy to see that $\mathrm{Pres}^{n}_{_{C}\mathcal{X}}(S)$ satisfies the assumption of Lemma \ref{ATWC}.
Then there is a triangle $U\longrightarrow V\longrightarrow C\longrightarrow$ with
$U\in\mathrm{Pres}^{n}_{_{C}\mathcal{X}}(S)$ and $V \in (\widecheck{\mathrm{Add}S})_{n}$.
Since $U\in \mathrm{Pres}^{n}_{_{C}\mathcal{X}}(S)=S^{\bot_{>0}}\bigcap$ $_{C}\mathcal{X}\subseteq C^{\bot_{>0}}$,
the triangle $U\longrightarrow V\longrightarrow C\longrightarrow$
is split, and then  $C\in (\widecheck{\mathrm{Add}S})_{n}$ by Lemma \ref{WEI} (4).

(3) We claim that $\mathrm{Hom}_{\mathcal{D}}(S,~T[i>n])=0$ for any $T\in$ $_{C}\mathcal{X}$.

In fact, for any $T\in$ $_{C}\mathcal{X}$, there are triangles
$T_{i}\longrightarrow C_{i-1}\longrightarrow T_{i-1}\longrightarrow$ with
$C_{i}\in \mathrm{Add}C$ for any $1\leq i\leq n$, where $T_{0}=T$.
Note that $(\widecheck{\mathrm{Add}S})_{n}$ is closed under direct summands by Lemma \ref{WEI} (3').
So we have that $C_{i} \in (\widecheck{\mathrm{Add}S})_{n}$ by (2). By Lemma \ref{ATWI},
there is a triangle $T\longrightarrow V\longrightarrow U\longrightarrow$
with $U\in (\widecheck{\mathrm{Add}S})_{n-1}$ and there
are triangles $V_{i}\longrightarrow Z_{i+1}\longrightarrow V_{i+1}\longrightarrow$ with $Z_{i}\in \mathrm{Add}S$
for any $0\leq i\leq n-1$, where $V_{0}=T_{n}$ and $V=V_{n}$. Note that $T_{n}\in$ $_{C}\mathcal{X}$.
It is easy to see that $V=V_{n}\in \mathrm{Pres}^{n}_{_{C}\mathcal{X}}(S)=S^{\bot_{>0}}\bigcap$ $_{C}\mathcal{X}$.
Applying the functor $\mathrm{Hom}_{\mathcal{D}}(S,-)$ to the triangle $T\longrightarrow V\longrightarrow U\longrightarrow$,
we have that $\mathrm{Hom}_{\mathcal{D}}(S,U[i])=\mathrm{Hom}_{\mathcal{D}}(S,T[i+1])$ by the dimension shifting.
Note that $\mathrm{Hom}_{\mathcal{D}}(S,U[i>n-1])=0$ by Lemma \ref{WEI} (3'). Consequently,
$\mathrm{Hom}_{\mathcal{D}}(S,~T[i>n])=0$.

Since $S\in \mathrm{Pres}^{n}_{_{C}\mathcal{X}}(S)=S^{\bot_{>0}}\bigcap$ $_{C}\mathcal{X}$,
by Lemma \ref{WEI} (3), we have that $S\in (\widehat{\mathrm{Add}C})_{n}$.
\ \hfill $\Box$

\mskip\
Comparing to the characterization of tilting modules \cite[Theorem~3.1]{BS}, we have the following conclusion.

\bg{Th}\label{char}
Assume that $C$ is semi-selforthognal.
Then $(C,~S)$ is an n-silting pair if and only if $\mathrm{Pres}^{n}_{_{C}\mathcal{X}}(S)=S^{\bot_{>0}}\bigcap$ $_{C}\mathcal{X}$.

\ed{Th}

\Pf. By Proposition \ref{three} and Proposition \ref{PRES}.
\ \hfill $\Box$

\bg{Cor}\label{Cor1}
Assume that $C$ is semi-selforthognal.
Then $(C,~S)$ is an n-silting pair if and only if $\mathrm{Pres}^{n}_{_{C}\mathcal{X}}(S)=S^{\bot_{>0}}$ and $\mathrm{Add}S\subseteq$ $_{C}\mathcal{X}$.
\ed{Cor}

\Pf. $\Rightarrow$ We have that $\mathrm{Pres}^{n}_{_{C}\mathcal{X}}(S)=S^{\bot_{>0}}\bigcap$ $_{C}\mathcal{X}$
by Theorem \ref{char}. So we only need to prove that $S^{\bot_{>0}}\subseteq$ $_{C}\mathcal{X}$.
By Proposition \ref{PRES}, we can obtain that $C\in (\widecheck{\mathrm{Add}S})_{n}$, and then $S^{\bot_{>0}}\subseteq C^{\bot_{>0}}$
by Lemma \ref{LIIL}. Consequently, $S^{\bot_{>0}}\subseteq$ $_{C}\mathcal{X}$ by Lemma 3.11 in \cite{WEI}.
By the lemma \ref{supplement}, $\mathrm{Add}S\subseteq$ $_{C}\mathcal{X}$.

$\Leftarrow$ We only need to prove that $S^{\bot_{>0}}\subseteq$ $_{C}\mathcal{X}$ by Theorem \ref{char}.
Since $S\in\mathrm{Pres}^{n}_{_{C}\mathcal{X}}(S)=S^{\bot_{>0}}$, $S$ is semi-selforthognal.
For any $X\in S^{\bot_{>0}}= \mathrm{Pres}^{n}_{_{C}\mathcal{X}}(S)$, then there are some triangles
$X_{i}\longrightarrow S_{i-1}\longrightarrow X_{i-1}\longrightarrow$ with $S_{i}\in \mathrm{Add}S$ for $1\leq i\leq n$
and $X_{n}\in$ $_{C}\mathcal{X}$, where $X_{0}=X$. Note that both $X_{n}$ and $\mathrm{Add}S$ are in $_{C}\mathcal{X}$.
Since $_{C}\mathcal{X}$ is closed under extensions and [1], $X_{n-1}\in$ $_{C}\mathcal{X}$.
Similar, we can get that$X\in$ $_{C}\mathcal{X}$. i.e., $S^{\bot_{>0}}\subseteq$ $_{C}\mathcal{X}$.
\ \hfill $\Box$

\mskip\
Recall that a complex $T\in\mathcal{D}^{\leq0}$ is said to be semi-$n$-tilting, if $T$ is semi-tilting and
$T\in (\widehat{\mathrm{Add}R})_{n}$, see \cite[Definition~3.10]{WEI}.
Since $_{R}\mathcal{X}=R^{\bot_{>0}}=\mathcal{D}^{\leq0}$,
then we have the following result , see Theorem 4.4 in \cite{WEI}.
%
%
%
%
%
%
%

\bg{Cor}\label{Cor2}
Assume that $T\in \mathcal{D}^{\leq0}$. Then the following are equivalent:

$(1)$ T is semi-n-tilting;

$(2)$ $\mathrm{Pres}^{n}_{\mathcal{D}^{\leq0}}(T)=T^{\bot_{>0}}$
\ed{Cor}

\Pf. Note that $_{R}\mathcal{X}=R^{\bot_{>0}}=\mathcal{D}^{\leq0}$. By the theorem \ref{char}, we only need to prove that
$T^{\bot_{>0}}\subseteq$ $_{R}\mathcal{X}=R^{\bot_{>0}}$. This is easy by Lemma \ref{LIIL}.

\section{Relative AR-correspondence and co-t-structure}

Auslander and Reiten showed that there is a one-to-one correspondence between
isomorphism classes of basic cotilting  modules and certain contravariantly
finite resolving subcategories \cite{AR}.
Extending this result, Buan \cite{Buan} showed that there is a one-to-one correspondence between basic cotilting
complexes and certain contravariantly finite subcategories of the bounded derived category
of an artin algebra. In \cite{WZ1}, the authors extended such a result to cosilting complexes. In this section,
the author's mainly proof consists in that there is also an AR-correspondence with respect to the silting pair.
Furthermore, the authors gives a bijection between silting pair and bounded above co-t-structure.

Firstly, we give the following notions. Let $\mathcal{A}\subseteq\mathcal{B}$ be two subcategories of $\mc{D}$.
$\mathcal{A}$ is said to be covariantly finite in $\mathcal{B}$,
if for any $B\in \mathcal{B}$, there is a homomorphism $f$ : $B\to A$ for some $A\in \mathcal{A}$ such that
$\mathrm{Hom}_{\mathcal{D}}(f,~A')$ is surjective for any $A'\in \mathcal{A}$. Moreover, $\mathcal{A}$ is said to be specially covariantly finite in $\mathcal{B}$,
if for any $B\in \mathcal{B}$, there is a triangle $B\to A\to C\to$ with some $A\in \mathcal{A}$ such that
$\mathrm{Hom}_{\mathcal{D}}(C,~A'[1])=0$ for any $A'\in\mathcal{A}$. Note that in the latter case, one has that
$C\in$ $^{\bot_{>0}}\mathcal{A}$ if $\mathcal{A}$ is closed under $[1]$.

Recall that $\mathcal{X}*\mathcal{Y}$ denotes the class of objects $Z \in \mathcal{D}$ for which there is a triangle
$X\longrightarrow Z\longrightarrow Y\longrightarrow$ with $X\in \mathcal{X}$ and $Y\in \mathcal{Y}$. A pair $(\mathcal{A},~\mathcal{B})$ of
subcategories in $\mathcal{D}$ is said to be a co-t-structure on $\mathcal{D}$ (see \cite{PD}) if the following statements hold: (1) $\mathcal{A}$ is closed under [-1] and
$\mathcal{B}$ is closed under [1]; (2) $\mathrm{Hom}_{\mathcal{D}}(\mathcal{A}[-1],~\mathcal{B})=0$; (3) $\mathcal{D}=\mathcal{A}[-1]*\mathcal{B}$.
Moreover, the co-t-structure $(\mathcal{A},~\mathcal{B})$ on $\mathcal{D}$ is said to be bounded above (see \cite{HVVS}) if $\bigcup_{i\in \mathrm{Z}}\mathcal{B}[i]=\mathcal{D}$.

\mskip\

The following lemma is the dual of Corollary 2.3 in \cite{WZ1}. For the reader's convenience, we give a proof here.

\bg{Lem}\label{JIEGUO1}
Assume that S is semi-selforthogonal. For any $L\in \widecheck{_{S}\mathcal{X}}$,
there are two triangles $X\longrightarrow Y\longrightarrow L\longrightarrow$ with
$X\in$ $_{S}\mathcal{X}$, $Y\in (\widecheck{\mathrm{Add}S})_{n}$ for some integer $n$, and
$L\longrightarrow U\longrightarrow V\longrightarrow$ with
$V\in (\widecheck{\mathrm{Add}S})_{n-1}$ and $U\in S^{\bot_{>0}}$.
\ed{Lem}

\Pf. Since $L\in \widecheck{_{S}\mathcal{X}}$, for some integer $n$, there exist some triangles
$L_{i}\longrightarrow M_{i}\longrightarrow L_{i+1}\longrightarrow$ with
$M_{i}\in$ $_{S}\mathcal{X}$ for any $0\leq i\leq n$, where $L_{0}=L$ and $L_{n+1}=0$.
Thus we have a triangle $X\longrightarrow Y\longrightarrow L\longrightarrow$ with
$X\in$ $_{S}\mathcal{X}$, $Y\in (\widecheck{\mathrm{Add}S})_{n}$ by Lemma 2.7 in \cite{WEI}.
And then there is a triangle $Y\longrightarrow B_{0}\longrightarrow V\longrightarrow$ with
$V\in (\widecheck{\mathrm{Add}S})_{n-1}$ and $B_{0}\in \mathrm{Add}S$.
We have the following triangles commutative diagram:

$$\xymatrix{
X\ar@^{=}[d]\ar[r]&Y\ar[d]\ar[r]&L\ar[d]\ar[r]&\\
X\ar[r]&B_{0}\ar[r]\ar[d]&U\ar[d]\ar[r]&\\
&V\ar@^{=}[r]\ar[d]&V\ar[d]&\\
&&
}$$
Note that $X\in$ $_{S}\mathcal{X} \subseteq S^{\bot_{>0}}$ and $B_{0}\in S^{\bot_{>0}}$.
So the first row and the third column are just the desired triangles.
\ \hfill $\Box$

\bg{Pro}\label{pro1}
If the pair $(C,~S)$ is silting, then we have the following three conclusions:

$(1)$ $\widecheck{S^{\bot_{>0}}}=C^{\bot_{\gg0}}$;

$(2)$ $S^{\bot_{>0}}$ is specially covariantly finite in $C^{\bot_{\gg0}}$;

$(3)$ the pair $(^{\bot_{0}}(S^{\bot_{>0}})[1]\bigcap C^{\bot_{\gg0}},S^{\bot_{>0}})$ is a bounded above co-t-structure on $C^{\bot_{\gg0}}$.

\ed{Pro}

\Pf. (1) For any $M\in\widecheck{S^{\bot_{>0}}}$, there are some triangles
$M_{i-1}\longrightarrow X_{i-1}\longrightarrow M_{i}\longrightarrow$ with
$X_{i}\in S^{\bot_{>0}}$ for $1\leq i\leq n+1$, where $M_{0}=M$ and $M_{n+1}=0$.
Note that $X_{i}\in S^{\bot_{>0}}\subseteq C^{\bot_{>0}}$ by Lemma \ref{LIIL}.
Applying the functor $\mathrm{Hom}_{\mathcal{D}}(C,~-)$
to these triangles, we have that
$0=\mathrm{Hom}_{\mathcal{D}}(C,~X_{n}[1])\cong \mathrm{Hom}_{\mathcal{D}}(C,~M_{n}[1])\cong
\cdots \cong\mathrm{Hom}_{\mathcal{D}}(C,~M[n+1])$ by dimension shifting,
and then $M\in C^{\bot_{\gg0}}$. i.e., $\widecheck{S^{\bot_{>0}}}\subseteq C^{\bot_{\gg0}}$.

For the reverse inclusion, if $N\in C^{\bot_{\gg0}}$, then we have that $L=N[t]\in C^{\bot_{>0}}$ for some $t>0$.
We claim that there exists some $s>0$ such that $L[s]\in S^{\bot_{>0}}$.
Indeed, since $S\in \widehat{\mathrm{Add}C}$, there exist some triangles
$S_{i}\longrightarrow C_{i-1}\longrightarrow S_{i-1}\longrightarrow$ with $C_{i}\in \mathrm{Add}C$
for any $1\leq i\leq s+1$, where $S_{0}=S$ and $S_{s+1}=0$.
Applying the functor $\mathrm{Hom}_{\mathcal{D}}(-,~L)$ to these triangles, we have that
$0=\mathrm{Hom}_{\mathcal{D}}(C_{s},~L[1])\cong \mathrm{Hom}_{\mathcal{D}}(S_{s},~L[1])\cong
\cdots \cong \mathrm{Hom}_{\mathcal{D}}(S,~L[s+1])$ by dimension shifting,
and then $L\in S^{\bot_{>s}}$. i.e., $L[s]\in S^{\bot_{>0}}$ for some $s>0$.
So $N[s+t]\in S^{\bot_{>0}}$. Since $\widecheck{S^{\bot_{>0}}}$ is closed under [-1] and
$S^{\bot_{>0}}\subseteq\widecheck{S^{\bot_{>0}}}$,
$N\in \widecheck{S^{\bot_{>0}}}$. i.e., $C^{\bot_{\gg0}}\subseteq\widecheck{S^{\bot_{>0}}}$.

(2) For any $P\in C^{\bot_{\gg0}}$, by Lemma 3.11 in \cite{WEI} and (1), we have that $P\in\widecheck{_{S}\mathcal{X}}$.
By the lemma \ref{JIEGUO1}, there is a triangle
$P\longrightarrow U\longrightarrow V\longrightarrow$ with
$V\in \widecheck{\mathrm{Add}S}$ and $U\in S^{\bot_{>0}}$.
We only need to prove that $\mathrm{Hom}_{\mathcal{D}}(V,~S^{\bot_{>0}}[1])=0$ by the definition.
Since $V\in \widecheck{\mathrm{Add}S}$, there exist some triangles
$V_{i}\longrightarrow S_{i}\longrightarrow V_{i+1}\longrightarrow$ with
$S_{i}\in \mathrm{Add}S$ for any $0\leq i\leq n$, where $V_{0}=V$ and $V_{n+1}=0$.
It is not difficult to verify that $\mathrm{Hom}_{\mathcal{D}}(V,~S^{\bot_{>0}}[1])=0$
by dimension shifting. Consequently,
$S^{\bot_{>0}}$ is specially covariantly finite in $C^{\bot_{\gg0}}$.

(3) Note that $C^{\bot_{\gg0}}$ is a triangulated subcategory of $\mathcal{D}$. It is easy to see that $^{\bot_{0}}(S^{\bot_{>0}})[1]\bigcap C^{\bot_{\gg0}}$ is
closed under [-1], $S^{\bot_{>0}}$ is closed under [1]
and $\mathrm{Hom}_{\mathcal{D}}((^{\bot_{0}}(S^{\bot_{>0}})[1]\bigcap C^{\bot_{\gg0}}[-1],~S^{\bot_{>0}})=0$.
For any $X\in C^{\bot_{\gg0}}$, there is a triangle $A[-1]\longrightarrow X\longrightarrow B\longrightarrow A$ with
$B\in S^{\bot_{>0}}$ and $\mathrm{Hom}_{\mathcal{D}}(A,~S^{\bot_{>0}}[1])=0$ by (2),
then $A\in$ $^{\bot_{0}}(S^{\bot_{>0}})[1]\bigcap C^{\bot_{\gg0}}$.
So $C^{\bot_{\gg0}}=$ $^{\bot_{0}}(S^{\bot_{>0}})[1]\bigcap C^{\bot_{\gg0}}*S^{\bot_{>0}}$,
i.e., $(^{\bot_{0}}(S^{\bot_{>0}})[1]\bigcap C^{\bot_{\gg0}},S^{\bot_{>0}})$ is a co-t-structure on $C^{\bot_{\gg0}}$.

It is not difficult to verify that $S^{\bot_{>0}}$ is suspended. By Lemma 2.1 in \cite{WEI1}, we have
$$C^{\bot_{\gg0}}=\widecheck{S^{\bot_{>0}}}=\bigcup_{i\in Z}S^{\bot_{>0}}[i].$$
Thus, it is bounded above.
\ \hfill $\Box$

\mskip\

\bg{Lem}\label{JIEGUO3}
If both $(C,~S)$ and $(C,~S\oplus T)$ are silting pairs, then $T\in \mathrm{Add}S$.

\ed{Lem}

\Pf. Since both  $S$ and $S\oplus T$ are semi-selforthogonal, $T$ is also semi-selforthogonal
and $\mathrm{Hom}_{\mathcal{D}}(S,~T^{(I)}[i>0])=0$ for any set $I$. Note that there is a triangle
$T_{1}\longrightarrow S^{(X_{0})}\stackrel{u_{0}}{\longrightarrow} T\longrightarrow$ with
$u_{0}$ evaluation map. It follows that $T_{1}\in S^{\bot_{>0}}$. Similarly, there is a triangle
$T_{2}\longrightarrow S^{(X_{1})}\stackrel{u_{1}}{\longrightarrow} T_{1}\longrightarrow$ with
$u_{1}$ evaluation map and $T_{2}\in S^{\bot_{>0}}$. Repeat the process, we have that $T\in$ $_{S}\mathcal{X}$.

Since both $S$ and $S\oplus T$ are in $\widehat{\mathrm{Add}C}$, $T$ is in $\widehat{\mathrm{Add}C}$.
Thus there exist some triangles $T_{i+1}\longrightarrow C_{i}\longrightarrow T_{i}\longrightarrow$
with $C_{i}\in\mathrm{Add}C$ for any $0\leq i\leq n$, where $T_{0}=T$ and $T_{n+1}=0$.
For any $M\in S^{\bot_{>0}}\subseteq C^{\bot_{>0}}$, we apply the functor $\mathrm{Hom}_{\mathcal{D}}(-,~M)$ to these triangles,
and obtain that $0=\mathrm{Hom}_{\mathcal{D}}(C_{n},~M[i])\cong \mathrm{Hom}_{\mathcal{D}}(T_{n},~M[i])\cong \cdots
\cong \mathrm{Hom}_{\mathcal{D}}(T,~M[n+i])$ for any $i>0$. i.e., $T\in$ $^{\bot_{>n}}(S^{\bot_{>0}})$.
It follows from Lemma \ref{WEI} (3) that $T\in (\widehat{\mathrm{Add}S})_{n}$.
Since $S\oplus T$ is semi-selforthogonal, $T\in$ $^{\bot_{>0}}(\mathrm{Add}S)$.
Thus $T\in \mathrm{Add}S$ by Lemma \ref{WEI} (5).
\ \hfill $\Box$

\bg{Pro}\label{pro2}
Assume that C is semi-selforthogonal, $\mathcal{Y}\subseteq C^{\bot_{>0}}$ is specially covariantly finite in $C^{\bot_{\gg0}}$
and is suspended such that $\widecheck{\mathcal{Y}}=C^{\bot_{\gg0}}$.
If $\mathcal{Y}\bigcap$ $^{\bot_{>0}}\mathcal{Y}$ is closed under direct sums,
then there is a silting pair $(C,~S)$ such that $\mathcal{Y}=S^{\bot_{>0}}$.

\ed{Pro}

\Pf. For any $X\in \widecheck{\mathcal{Y}}$, set $\mathcal{Y}$-cores.dim $(X)=s$. i.e., there are some triangles
$X_{i}\longrightarrow Y_{i}\longrightarrow X_{i+1}\longrightarrow$
with $Y_{i}\in \mathcal{Y}$ for any $0\leq i\leq s-1$, where $X_{0}=X$.
By the dimension shifting, it is not difficult to verify that
$X\in (^{\bot_{>0}}\mathcal{Y})^{\bot_{>s}}$.

Since $\mathcal{Y}\subseteq C^{\bot_{>0}}$ is specially covariantly finite in $C^{\bot_{\gg0}}$
and is suspended, for any $M\in C^{\bot_{\gg0}}=\widecheck{\mathcal{Y}}$,
there are some triangles $M_{i}\longrightarrow Y_{i}\longrightarrow M_{i+1}\longrightarrow$
with $Y_{i} \in \mathcal{Y}$ and $M_{j}\in$ $^{\bot_{>0}}\mathcal{Y}$
(Look at the last sentence of the second paragraph at the beginning of the third section) for any $i\geq0$ and $j>0$,
where $M_{0}=M$. We claim that $M_{n}\in$ $\mathcal{Y}\bigcap$ $^{\bot_{>0}}\mathcal{Y}$, where $n=$ $\mathcal{Y}$-cores.dim $(M)$.
In fact, applying the functor $\mathrm{Hom}_{\mathcal{D}}(M_{n+1},~-)$ to these
triangles, we have that $\mathrm{Hom}_{\mathcal{D}}(M_{n+1},~M[i])\cong \mathrm{Hom}_{\mathcal{D}}(M_{n+1},~M_{n}[i-n])$
for any $i>n$ by the dimension shifting. Since $M_{n+1}\in$ $^{\bot_{>0}}\mathcal{Y}$ and $M\in (^{\bot_{>0}}\mathcal{Y})^{\bot_{>n}}$,
$\mathrm{Hom}_{\mathcal{D}}(M_{n+1},~M_{n}[1])=0$. i.e., the triangle $M_{n}\longrightarrow Y_{n}\longrightarrow M_{n+1}\longrightarrow$
is split. Since $\mathcal{Y}\bigcap$ $^{\bot_{>0}}\mathcal{Y}$ is closed under direct sums and $\mathcal{Y}$
is suspended, $\mathcal{Y}\bigcap$ $^{\bot_{>0}}\mathcal{Y}$ is closed under direct summands.
In fact, set $N=N_{1}\bigoplus N_{2}\in\mathcal{Y}\bigcap$ $^{\bot_{>0}}\mathcal{Y}$,
then $N^{(X)}=N_{1}\bigoplus N_{2}\bigoplus N_{1}\bigoplus N_{2}\cdots=N_{1}\bigoplus N^{(X)}$ for any set $X$.
Thus there is a split triangle $N^{(X)}\longrightarrow N^{(X)}\longrightarrow N_{1}\longrightarrow$.
Since $\mathcal{Y}$ is suspended, $N_{1}\in\mathcal{Y}$, and then $N_{1}\in\mathcal{Y}\bigcap$ $^{\bot_{>0}}\mathcal{Y}$.
i.e., $\mathcal{Y}\bigcap$ $^{\bot_{>0}}\mathcal{Y}$ is closed under direct summands.
So $M_{n}\in$ $\mathcal{Y}\bigcap$ $^{\bot_{>0}}\mathcal{Y}$.

Note that $C\in C^{\bot_{\gg0}}=\widecheck{\mathcal{Y}}$. Specially, we take $M=C$ in the discussion of the second paragraph.
There are triangles $C_{i}\longrightarrow Y_{i}\longrightarrow C_{i+1}\longrightarrow$
with $Y_{i} \in \mathcal{Y}\bigcap$ $^{\bot_{>0}}\mathcal{Y}$ and $C_{i}\in$ $^{\bot_{>0}}\mathcal{Y}$ for any $0\leq i\leq m$,
where $C_{0}=C$ and $C_{m+1}=0$. Note that $\mathcal{Y}\subseteq C^{\bot_{>0}}$. i.e., $C\in$ $^{\bot_{>0}}\mathcal{Y}$.
Set $S=\bigoplus_{j=0}^{m}Y_{j}$. Then $S$ is semi-selforthogonal since
$\mathcal{Y}\bigcap$ $^{\bot_{>0}}\mathcal{Y}$ is closed under direct sums.
We claim that $(C,~S)$ is a silting pair. We only need to prove that $S\in \widehat{\mathrm{Add}C}$ by the definition.
There are some triangles $S_{i+1}\longrightarrow C^{(X_{i})}\stackrel{u_{i}}{\longrightarrow} S_{i}\longrightarrow$ with
$u_{i}$ evaluation map for any $i\geq0$, where $S_{0}=S$. Since $S\in \mathcal{Y}\subseteq C^{\bot_{>0}}$ and
$C$ is semi-selforthogonal, $S_{i} \in C^{\bot_{>0}}$ for any $i\geq1$. Thus $S\in$ $_{C}\mathcal{X}$.
If $S\in$ $^{\bot_{>t}}(C^{\bot_{>0}})$ for some $t$, then $S\in \widehat{\mathrm{Add}C}$ by Lemma \ref{WEI} (3).
Indeed, for any $T\in C^{\bot_{>0}}\subseteq C^{\bot_{\gg0}}=\widecheck{\mathcal{Y}}$, there are some triangles
$T_{i}\longrightarrow Z_{i}\longrightarrow T_{i+1}\longrightarrow$
with $Z_{i} \in \mathcal{Y}$ for any $0\leq i\leq t$, where $T_{0}=T$ and $T_{t+1}=0$.
Applying the functor $\mathrm{Hom}_{\mathcal{D}}(S,~-)$ to these triangles, we have that
$\mathrm{Hom}_{\mathcal{D}}(S,~T[i])\cong \mathrm{Hom}_{\mathcal{D}}(S,~Z_{t}[i-t])$ for any $i>t$ by the dimension
shifting. Since $S=\bigoplus_{j=0}^{m}Y_{j}\in$ $^{\bot_{>0}}\mathcal{Y}$ and $Z_{t}\in \mathcal{Y}$,
$\mathrm{Hom}_{\mathcal{D}}(S,~T[i>t])=0$. i.e., $S\in$ $^{\bot_{>t}}(C^{\bot_{>0}})$.
Consequently, $(C,~S)$ is a silting pair.

Finally, we prove that $\mathcal{Y}=S^{\bot_{>0}}$.
Since $S=\bigoplus_{j=0}^{m}Y_{j}\in$ $^{\bot_{>0}}\mathcal{Y}$, we have that $\mathcal{Y}\subseteq S^{\bot_{>0}}$.
For reverse inclusion, for any $T\in S^{\bot_{>0}}\subseteq C^{\bot_{>0}}$ by Lemma \ref{LIIL} and by the discussion above,
there exist some triangles $T_{i}\longrightarrow X_{i}\longrightarrow T_{i+1}\longrightarrow$
with $X_{i} \in \mathcal{Y}\bigcap$ $^{\bot_{>0}}\mathcal{Y}$, for any $0\leq i\leq n$,
where $T_{0}=T$ and $T_{n+1}=0$. Note that $\mathcal{Y}\subseteq S^{\bot_{>0}}$. Since $S^{\bot_{>0}}$ is
suspended, all objects in these triangles are in $S^{\bot_{>0}}$. It is not difficult to verify that
$(C,~S\oplus U)$ is a silting pair for any $U \in\mathcal{Y}\bigcap$ $^{\bot_{>0}}\mathcal{Y}$.
Indeed, we only need to prove that $S\oplus U\in \widehat{\mathrm{Add}C}$. This prove similar to $S\in \widehat{\mathrm{Add}C}$.
Thus $U\in\mathrm{Add}S$ by Lemma \ref{JIEGUO3}, and then $\mathcal{Y}\bigcap$ $^{\bot_{>0}}\mathcal{Y}$ $\subseteq \mathrm{Add}S$.
Note that $\mathcal{Y}\bigcap$ $^{\bot_{>0}}\mathcal{Y}$ is closed under direct sums and direct summands.
It is easy to see that $\mathcal{Y}\bigcap$ $^{\bot_{>0}}\mathcal{Y}=\mathrm{Add}S$.
So $T_{1}\in\widecheck{\mathrm{Add}S}\bigcap S^{\bot_{>0}}=\mathrm{Add}S$.
Then the triangle $T\longrightarrow X_{0}\longrightarrow T_{1}\longrightarrow$ is split.
Since $T_{1}$, $X_{0}$ $\in \mathcal{Y}$ and $\mathcal{Y}$ is suspended, we have that $T\in\mathcal{Y}$.
So $\mathcal{Y}=S^{\bot_{>0}}$.

Thus we complete this proof.
\ \hfill $\Box$

\bg{Cor}\label{cot2}
Assume that C is semi-selforthogonal and the subcategory $\mathcal{Y}\subseteq C^{\bot_{>0}}$ is closed under direct sums.
If $(\mathcal{X},~\mathcal{Y})$ is a bounded above co-t-structure on $C^{\bot_{\gg0}}$, then there is a silting pair $(C,~S)$ such that
$\mathcal{Y}=S^{\bot_{>0}}$.

\ed{Cor}

\Pf. For any $A\in C^{\bot_{\gg0}}$, there is a triangle $X[-1]\longrightarrow A\longrightarrow Y\longrightarrow X$
with $X\in \mathcal{X}$ and $Y\in \mathcal{Y}$. Obviously, $\mathrm{Hom}_{\mathcal{D}}(X,~Y'[1])=0$ for all $Y'\in\mathcal{Y}$. Thus $\mathcal{Y}$ is
specially covariantly finite in $C^{\bot_{\gg0}}$. It follows from Proposition 2.1 in \cite{PD} that $\mathcal{Y}$ is
suspended. By the lemma 2.1 in \cite{WEI1} and the assumption, we have
$$\widecheck{\mathcal{Y}}=\bigcup_{i\in Z}\mathcal{Y}[i]=C^{\bot_{\gg0}}.$$
Consequently, $\mathcal{Y}$ satisfies the conditions of Proposition \ref{pro2}, and then
there is a silting pair $(C,~S)$ such that $\mathcal{Y}=S^{\bot_{>0}}$.
\ \hfill $\Box$

\bg{Lem}\label{JIEGUO4}
If the pair $(C,~S)$ is an n-silting pair, then $S^{\bot_{>0}}=\mathrm{Pres}^{n}_{C^{\bot_{>0}}}(S)$.

\ed{Lem}

\Pf. By the lemma 3.11 from \cite{WEI} and Proposition \ref{three}.
%
%
%
\ \hfill $\Box$

\mskip\

Two $n$-silting pairs $(C,~S)$ and $(C,~T)$ are said to be equivalent if $\mathrm{Add}S=\mathrm{Add}T$.
By Proposition \ref{pro1} and Proposition \ref{pro2}, we can obtain the AR-correspondence with
respect to silting pair as follows.

\bg{Th}\label{ARC}
If $C$ is semi-selforthogonal, then
there is an one-to-one correspondence between equivalent class of n-silting pairs $(C,~S)$ and subcategories
$\mathcal{Y}\subseteq C^{\bot_{>0}}$ which is specially covariantly finite in $C^{\bot_{\gg0}}$,
suspended and closed under direct sums such that $\widecheck{\mathcal{Y}}=C^{\bot_{\gg0}}$,
given by $F$: $(C,~S)$ $\longrightarrow S^{\bot_{>0}}$.
\ed{Th}

\Pf. It follows from Proposition \ref{pro1}, Proposition \ref{pro2} and Lemma \ref{JIEGUO4} that the correspondence is well-defined.
Indeed, if the two $n$-silting pairs $(C,~S)$ and $(C,~T)$ are equivalent, by Lemma \ref{JIEGUO4}, we have that
$S^{\bot_{>0}}=\mathrm{Pres}^{n}_{C^{\bot_{>0}}}(S)=\mathrm{Pres}^{n}_{C^{\bot_{>0}}}(T)=T^{\bot_{>0}}$.
Moreover, $F$ is surjective by Proposition \ref{pro2}.
If both $(C,~S)$ and $(C,~T)$ are $n$-silting pairs with $S^{\bot_{>0}}=T^{\bot_{>0}}$,
it is easy to verify that $(C,~S\oplus T)$ is also $n$-silting pair by the definition.
So we have that $\mathrm{Add}S=\mathrm{Add}T$ by Lemma \ref{JIEGUO3}, i.e., $(C,~S)$ and $(C,~T)$ are equivalent.
Hence, $F$ is bijective.
\ \hfill$\Box$

\mskip\

In \cite{HVVS}, the authors give a bijective correspondence between bounded co-t-structure and a silting class
(i.e., semi-selforthogonal class), see Theorem 5.10.
From the Proposition \ref{pro1} (3), Lemma \ref{JIEGUO4} and Corollary \ref{cot2}, we can obtain the following result.
It is regrettable that the co-t-structure in the following result is bounded above.

\bg{Th}\label{cot3}
If $C$ is semi-selforthogonal, then
there is a bijection between equivalent class of silting pairs $(C,~S)$ and bounded
above co-t-structure $(\mathcal{X},~\mathcal{Y})$ on $C^{\bot_{\gg0}}$ with
$\mathcal{Y}\subseteq C^{\bot_{>0}}$ closed under direct sums.

\ed{Th}

\Pf. This proof is similar to the proof of Theorem \ref{ARC}.
\ \hfill$\Box$

\section{Acknowledgments}

The authors would like to thank the referee for his/her careful reading of the paper and useful suggestions.

{\small

}

\end{document}